\input epsf
\font\sc=cmcsc10
\font\tencmmib=cmmib10
\font\eightcmmib=cmmib10
\font\tenmsym=msbm10    
\font\eightmsym=msbm8   
\font\teneufm=eufm10
\font\tencmsy=cmsy10

\font\twelvecmsy=cmmib10 scaled\magstep2
\textfont 7=\twelvecmsy \scriptfont 7=\tencmsy

\textfont 9=\tenmsym \scriptfont 9=\eightmsym
\def\bb {\fam9 }
\textfont 8=\teneufm \scriptfont 8=\teneufm
\def\frak {\fam8 }
\textfont 6=\tencmmib \scriptfont 6=\eightcmmib
\def\bi {\fam6 }
\magnification=1200
\font\bgf=cmbx10 scaled\magstep2

$\ $\vskip0.8in
\overfullrule=0pt

\centerline{\bgf Singular deformations of Lie algebras}\smallskip

\centerline{\bgf on an example}\bigskip

\centerline{\sc Alice FIALOWSKI \qquad and \qquad Dmitry FUCHS}\medskip
 
\noindent
\hbox{\hskip50pt
\vtop{
   \hbox to 2in {\hfil Department of Applied Analysis\hfil}\par
   \hbox to 2in {\hfil E\"otv\"os Lor\'and University\hfil} \par
   \hbox to 2in {\hfil M\'uzeum krt.~6-8\hfil}\par
   \hbox to 2in {\hfil H-1088 Budapest\hfil} \par
   \hbox to 2in {\hfil Hungary\hfil}}
\quad
\vtop{
   \hbox to 2in{\hfil Department of Mathematics\hfil} \par
   \hbox to 2in {\hfil University of California\hfil} \par
   \hbox to 2in {\hfil Davis CA 95616\hfil} \par
   \hbox to 2in {\hfil USA\hfil}}
}
 
\vskip 25pt
\centerline{\bf 1. Introduction}\medskip

We introduce the concept of singular deformations of Lie algebras.  As
far as we know they  have been never considered in the literature.
Nevertheless, they are unavoidable in any complete classification of
deformations.  In this paper we will show that singular deformations
occur even for some of the simplest infinite dimensional cases.

Let $\frak g$ be a Lie algebra  with the commutator $[\, ,\, ]$. Consider
a formal one-parameter deformation
 $$
[g,h]_t=[g,h]+\sum_{k\ge1}\alpha_k(g,h)t^k
 $$
of $\frak g$.  A deformation is called {\it non-singular} if there exists
a formal one-parameter family of linear transformations
 $$
\varphi_t(g)=g+\sum_{l\ge1}\beta_l(g)t^l
 $$
of $\frak g$ and a formal (not necessarily invertible) parameter change
$u=u(t)$ which transform the deformation $[g,h]_t$ into a deformation
 $$
[g,h]'_u=[g,h]+\sum_{k\ge1}\alpha'_k(g,h)u^k,\qquad
\varphi_t^{-1}[\varphi_t(g),\varphi_t(h)]_t=[g,h]'_{u(t)}
 $$
with the cocycle $\alpha'_1\in C^2({\frak g};{\frak g})$ being not
cohomologous to 0. Otherwise the deformation is called {\it singular}.

The example we present is the following.  Consider the complex
infinite-dimensional Lie algebra $L_1$ of polynomial vector fields in
$\bb C$ with trivial 1-jet at 0. This Lie algebra is spanned by the
vector fields $e_i=z^{i+1}\displaystyle{d\strut\over\strut dz},\,
i=1,2,3,\dots$, and the commutator is defined by the standard formula
$$[e_i,e_j]=(j-i)e_{i+j}.$$ This Lie algebra proves to be especially
interesting from the point of view of the deformation theory, for on one
hand its deformations may be completely classified, and on the other hand
they behave in a very unusual manner.

The deformations of $L_1$ were first studied in 1983 by the first author
[Fi1]. In [Fi1] three one-parameter deformations of the Lie algebra $L_1$
were considered:
 $$
\eqalign{[e_i,e_j]^1_t&=(j-i)(e_{i+j}+te_{i+j-1});\cr
[e_i,e_j]^2_t&=\cases{(j-i)e_{i+j}&if $i\ne1,j\ne1$,\cr (j-1)e_{j+1}+tje_j&if
$i=1,j\ne1;$\cr}\cr [e_i,e_j]^3_t&=\cases{(j-i)e_{i+j}&if $i\ne2,j\ne2$,\cr
(j-2)e_{j+2}+tje_j&if $i=2,j\ne2$.\cr}\cr}
 $$
All the three  families of Lie algebras may be realized as families of
subalgebras of the Lie algebra $L_0$ (spanned by $e_i$ with $i\ge0$). The
first deformation may be defined by the formula $e_i\mapsto e_i+te_{i-
1}\, (i\ge1)$; in other words, the Lie algebra $L_1$ of vector fields
with a double zero at 0 is deformed into the Lie algebra of vector fields
with two zeroes at points 0 and $t$. The two other deformations are
defined by the formulas $$\eqalign{&e_1\mapsto e_1+te_0,\ e_i\mapsto e_i\
{\rm if}\ i\ne1;\cr &e_2\mapsto e_2+te_0,\ e_i\mapsto e_i\ {\rm if}\
i\ne2.\cr}$$ All the three deformations are pairwise not equivalent.
Moreover, if $L_1^1,L_1^2,L_1^3$ are Lie algebras from the three families
corresponding to arbitrary non-zero values of the parameter (up to an
isomorphism, they do not depend on the non-zero parameter value), then
neither two of $L_1^1,L_1^2,L_1^3$ are isomorphic to each other. Indeed,
obviously $${\rm dim}\, (L_1^r/[L_1^r,L_1^r])=\cases{2&if $r=1$,\cr 1&if
$r=2,3$.\cr}$$ On the other hand, the Lie algebra $M^2=[L_1^2,L_1^2]$ is
spanned by $e_2,e_3,e_4,\dots$, and $M^3=[L_1^3,L_1^3]$ is spanned by
$e_1,e_3,e_4,\dots$. It is seen from this that
 $$
{\rm  dim}(M^r/[M^r,M^r])=\cases{3&if $r=2$,\cr 2&if $r=3$.\cr}
 $$

The main result of [Fi1] is the following\smallskip

{\sc Theorem 1.1}. {\it Any formal one-parameter deformation of $L_1$ may be
reduced by a formal parameter change to one of the three deformations
above.}\smallskip

However, the article [Fi1] contains no detailed proof of this
result. It is claimed there that the result follows from certain
calculations of Lie and Massey products in the cohomology
$H^\ast(L_1;L_1)$. These calculations are done in a more detailed
paper [Fi2], but still they do not imply Theorem 1.1. Moreover, the
description of the miniversal deformation of $L_1$ given in [Fi1]
needs a correction. Namely, the second of these three deformations is
{\it singular} in the above sense.  The correct description is the
following: the base of the miniversal deformation of $L_1$ is the
union of two smooth curves and one cuspidal curve passing through 0
with the common tangent.

\vbox{
\hbox to\hsize{\hfil\epsfbox{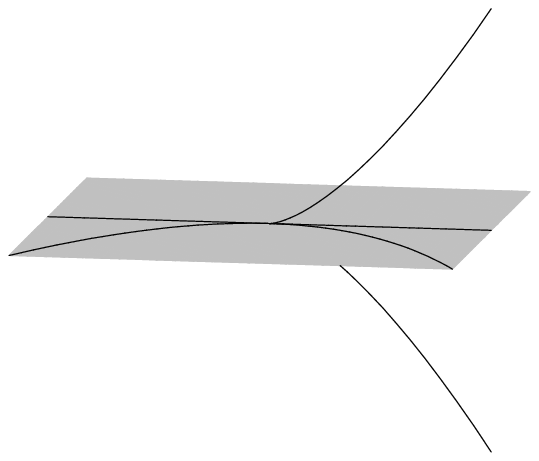}\hfil}
\centerline{Fig. 1}
}
\bigskip

Hence not only the base of the miniversal deformation of the Lie algebra
$L_1$ is a singular variety, but also one of its irreducible components
is singular. Although it is not impossible from the point of view of the
general theory, such examples have not been known before. As we have 
already mentioned, singular one parameter deformations, which appear, 
certainly, not only in the Lie algebra theory, have never been properly 
studied. For example, it is a common belief that if no non-trivial 
infinitesimal one parameter deformation is extendable to a formal 
deformation, then the Lie algebra is formally rigid. But is it 
really?\footnote*{Known (at least in the associative case) examples show 
that the relations between local (formal) and global (smooth) 
deformations may also be complicated in the infinite dimensional case; 
for example, the algebra of regular functions on the complex projective 
line with 4 punctures has a natural deformation with the cross ratio of 
the punctures as the parameter, but it has no infinitesimal deformations 
at all (see [Ko]).} 

There exists a general theory which provides a construction of a local
miniversal deformation of a Lie algebra $\frak g$. This theory is outlined in
[Fi3]. The article [Fi3] contains also a conjecture that the base of the
miniversal deformation may be described explicitly by a system of algebraic
equations in the cohomology space $H^2({\frak g};{\frak g})$, and that these
equations may be derived in a certain way from the Lie-Massey multiplication
in this cohomology. This procedure, however, needs some further
explanations. We are working on an explicit  general construction of
formal miniversal deformations of Lie algebras (see [FiFu]).

In this article we give direct proofs of Theorem 1.1 and of the above
description of the base of miniversal deformation of $L_1$. Both proofs are
direct and independent of any general theory; they appear completely 
reliable. 

In Sections 2 and 3 we list results which we regard as known; they concern
general theory of deformations of Lie algebras, and the cohomology of the Lie
algebra $L_1$. Section 4 contains the proof of the singularity of the 
deformation $[\, ,\, ]^2_t$ and the non-singularity of the deformations 
$[\, ,\, ]^1_t,[\, ,\, ]^3_t$. Theorem 1.1 is proved in Section 6 with 
some technical work done in Section 5. \bigskip 

\centerline{\bf 2. Deformations and cohomology}\medskip

Here we recall very briefly the classical theory of deformations of a Lie
algebra structure (see [Fu] and [Fi3] for details). Let $\frak g$ be a
(complex) Lie algebra with the bracket $[\, ,\, ]$. A formal 
one-parameter deformation of $\frak g$ is a power series 
$$[g,h]_t=[g,h]+\sum_{k\ge1}\alpha_k(g,h)t^k,$$ where $\alpha_k\in{\rm 
Hom}_{\bb C}(\Lambda^2{\frak g},{\frak g})=C^2({\frak g};{\frak g})$, the 
latter refers to the standard cochain complex $\{C^.({\frak g};{\frak 
g}),\, \delta\}$ of $\frak g$ with coefficients in the $\frak g$-module 
$\frak g$ (that is in the adjoint representation) -- see [Fu]. The Jacobi 
identity for $[\, ,\, ]_t$ is equivalent to the sequence of relations 
 $$
\delta\alpha_k=-{1\over2}\sum_{i=1}^{k-1}[\alpha_i,\alpha_{k-i}],\eqno(1) 
 $$ 
where $[\, ,\, ]$ denotes the usual product in the standard graded 
differential Lie algebra structure of the complex $\{C^.({\frak g};{\frak 
g}),\, \delta\}$: if $\beta\in C^p({\frak g};{\frak g})$ and $\gamma\in 
C^q({\frak g}; {\frak g})$ then $[\beta ,\gamma]\in C^{p+q-1}({\frak g}; 
{\frak g}),$ 
 $$
\eqalign{[\beta ,\gamma]&(g_1,\dots ,g_{p+q-1}) \cr &=
\sum_{1\le i_1<\dots <i_q\le p+q-1}(-1)^{i_1+\dots + i_q -
{q(q+1)\over2}}\beta(\gamma(g_{i_1},\dots ,g_{i_q}),g_1,\dots
 \hat g_{i_1}\dots\hat g_{i_q}\dots ,g_{p+q-1}) \cr
 &+(-1)^{pq+p+q}\sum_{1\le j_1<\dots <j_p\le p+q-1}(-1)^{j_1+\dots +j_p -
   {p(p+1)\over2}}\gamma(\beta(g_{j_1},\dots ,g_{j_p}),\cr 
 &\hskip3in g_1,\dots\hat g_{j_1}\dots\hat g_{j_p}\dots ,g_{p+q-
1}).\cr}
 $$
In particular,
 $$
\delta\alpha_1=0,\ \delta\alpha_2=-{1\over2}[\alpha_1,\alpha_1],\ 
\delta\alpha_3=-[\alpha_1,\alpha_2].\eqno(2) 
 $$ 
Two deformations $[g,h]_t,[g,h]'_t$ are 
called equivalent if there exists a formal one-parameter family 
$\{\varphi_t\}$ of linear transformations of $\frak 
g$,$$\varphi_t(g)=g+\sum_{k\ge1}\beta_k(g)t^k,$$such that 
$$[g,h]'_t=\varphi_t^{-1}[\varphi_t(g),\varphi_t(h)]_t.$$It is easy to 
see that $$\alpha'_1-\alpha_1=\delta\beta_1\eqno(3)$$and, more generally, 
if $\beta_1=\dots =\beta_{s-1}=0$, then$$\alpha'_s-
\alpha_s=\delta\beta_s$$ (here $\alpha'_i$ corresponds to $\alpha_i$ for 
$[\, ,\, ]'_t$). The first of equalities (2) shows that $\alpha_1$ is a 
cocycle, and the equality (3) shows that the cohomology class of 
$\alpha_1$ depends only on the equivalence class of the deformation. 
Cohomology classes from $H^2({\frak g};{\frak g})$ are called {\it 
infinitesimal deformations of} $\frak g$, and the cohomology class of 
$\alpha_1$ is called the {\it differential} of the formal deformation 
$[\, ,\, ]_t$. An infinitesimal deformation $\alpha\in H^2({\frak 
g};{\frak g})$ is not necessarily the differential of any formal 
deformation: the second and the third of the equalities (2) provide a 
necessary condition for it: the Lie square and the Massey-Lie cube should 
be equal to 0. The other relations (1) give more necessary conditions for 
an infinitesimal deformation being a differential (which comprise all 
together a sufficient condition for this). Usually these conditions are 
formulated in terms of higher Massey-Lie products (see [R], [Fu] and 
[FL]), but we use another method in this article.\bigskip 

\centerline{\bf 3. Cohomology of ${\bi L}_{\bf 1}$}\medskip

Here we recall the necessary information from [FeFu] and [Fi1] about the
(continuous) cohomology $H^\ast(L_1;L_1)$. First of all, the Lie algebra $L_1$
is ${\bb Z}_{>0}$-graded, ${\rm deg}\, e_i=i$, and this gives rise to a $\bb
Z$-grading in $C^\ast(L_1;L_1)$ and $H^\ast(L_1;L_1)$: ${\rm deg}\, \alpha=k$
for $\alpha\in C^q(L_1;L_1)$ if $$\alpha(e_{i_1},\dots ,e_{i_q})\in{\bb
C}e_{i_1+\dots +i_q-k}\ {\rm for\ all}\ i_1,\dots ,i_q.$$ One has
$H^q(L_1;L_1)=\oplus_kH^q_{(k)}(L_1;L_1)$. The following is a corollary from a
more general result of [FeFu] (see [Fi1]):\smallskip

{\sc Theorem 3.1}. {\it $$H^q_{(k)}(L_1;L_1)\cong\cases{\bb C&if
$\displaystyle{{3(q-1)^2+(q-1)\strut\over\strut2}\le k<{3q^2-q\over2}}$,\cr
0&otherwise.\cr}$$ In particular,
$$H^2(L_1;L_1)\cong\bigoplus_{k=2}^4H^2_{(k)}(L_1;L_1),\
H^3(L_1;L_1)\cong\bigoplus_{k=7}^{11}H^3_{(k)}(L_1;L_1),$$all the summands in
these two sums being isomorphic to $\bb C$.}\smallskip

Lie and Massey-Lie -products in $H^\ast(L_1;L_1)$ have been calculated in
[Fi2]. We need the following result from [Fi2].\smallskip

{\sc Theorem 3.2}. {\it Let $0\ne b\in H^2_{(3)}(L_1;L_1),\ 0\ne c\in
H^2_{(4)}(L_1;L_1)$. Then
 $$
\eqalign{&0\ne[b,c]\in H^3_{(7)}(L_1;L_1),\cr
         &0\ne[c,c]\in H^3_{(8)}(L_1;L_1),\cr
         &0\ne\langle b,b,b\rangle\in H^3_{(9)}(L_1;L_1).\cr}
 $$
}\smallskip

Here $\langle b,b,b\rangle$ is the Massey-Lie cube of $b$; the inequality
$\langle b,b,b\rangle\ne0$ means that if $\beta\in b$ is a cocycle, and
$[\beta,\beta]=\delta f$, then the cocycle $[\beta,f]$ is not
cohomologous to 0.

See the references for the proofs of Theorems 3.1 and 3.2. Theorem 3.2 may 
also be obtained from the results of Section 5 below.

\bigskip

\centerline{\bf 4. The three deformations of ${\bi L}_{\bf 1}$}\medskip

Consider the three deformations of the Lie algebra $L_1$ defined in the
introduction.\smallskip

{\sc Theorem 4.1.} {\it Of the three deformations $[\, ,\, ]_t^r,\, r=1,2,3,$
of $L_1$, the first and the third are non-singular, while the second is
singular.}\smallskip

{\sc Proof.} The three deformations have the form
$$[g,h]_t^r=[g,h]+\alpha_1^r(g,h)t,\,\qquad
r=1,2,3,$$where $$\eqalign{\alpha_1^1(e_i,e_j)&=(j-i)e_{i+j-1};\cr 
\alpha_1^2(e_i,e_j)&=\cases{je_j&if $i=1,j\ne1,$\cr 0&if 
$i\ne1,j\ne1;$\cr}\cr \alpha_1^3(e_i,e_j)&=\cases{je_j&if $i=2,j\ne2,$\cr 
0&if $i\ne2,j\ne2.$\cr}\cr}$$ Obviously, $\alpha_1^1,\alpha_1^2\in 
C^2_{(1)}(L_1;L_1),\, \alpha_1^3\in C^2_{(2)}(L_1;L_1)$. The cocycle 
$\alpha_1^3$ is not cohomologous to 0. To prove this, we calculate its 
value on a non-trivial cycle in $C^2_{(2)}(L_1;L_1^\ast),\ L_1^\ast$ 
being the $L_1$-module dual to $L_1$. This cycle may be chosen as 
$$a^{(2)}=e_3^\ast\otimes(e_1\wedge e_4-3e_2\wedge 
e_3)+{1\over2}e_2^\ast\otimes(e_1\wedge e_3)+3e_1^\ast\otimes(e_1\wedge 
e_2)$$ (the fact that it is a cycle is checked by a direct calculation, 
the fact that it is not homologous to 0 follows from the calculation 
below). We have: $$\alpha_1^3(e_1,e_4)=0,\, \alpha_1^3(e_2,e_3)=3e_3,\, 
\alpha_1^3(e_1,e_3)=0,\,  \alpha_1^3(e_1,e_2)=-
e_1,$$$$\langle\alpha_1^3,a^{(2)}\rangle=0-9+0-3=-12\ne0.$$ Hence the 
deformation $[\, ,\, ]^3_t$ is non-singular. 

The cocycles $\alpha_1^1,\alpha_1^2$ are cohomologous to 0 (because
$H^2_{(1)}(L_1;L_1)=0$). More precisely, $\alpha_1^r=\delta\beta^r,\, r=1,2,$
where $\beta^1,\beta^2\in C^1_{(1)}(L_1;L_1)$ are given by the formulas
$$\beta^1(e_i)={i-1\over2}e_{i-1},\
\beta^2(e_i)=\cases{\displaystyle{i+1\over2}e_{i-1}&for $i\ne1$,\cr 0&for
$i=1$.\cr}$$Put $\varphi^r_t(e_i)=e_i+\beta^r(e_i)t$, and
compute$$\gamma^r(e_i,e_j)=(\varphi^r_t)^{-1}[\varphi^r_t(e_i),\varphi^r_t
(e_j)]_t$$ modulo $t^4$. It follows from $\alpha_1^r=\delta\beta^r$ that
$$\gamma^r(e_i,e_j)=[e_i,e_j]+\sum_{k\ge2}\gamma^r_k(e_i,e_j)t^k,$$ where
$\gamma^r_k\in C^2_{(k)}(L_1;L_1)$, and a direct calculation shows that
$$\matrix{\gamma^1_2(e _1,e_2)=0,\hfill&\gamma^1_2(e_1,e_3)=e_2,\hfill&
\gamma^1_2(e_1,e_4)=3e_3,\hfill&\gamma^1_2(e_2,e_3)=\displaystyle{3\over2}e_3,
\cr \gamma^1_3(e_1,e_2)=0,\hfill&\gamma^1_3(e_1,e_3)=e_1,\hfill&
\gamma^1_3(e_1,e_4)=-3e_2,\hfill& \gamma^1_3(e_2,e_3)=-e_2,\cr
\gamma^2_2(e_1,e_2)=0,\hfill& \gamma^2_2(e_1,e_3)=4e_2,\hfill&
\gamma^2_2(e_1,e_4)=\displaystyle{15\over2}e_3,\hfill&
\gamma^2_2(e_2,e_3)=\displaystyle{15\over2}e_3,\cr
\gamma^2_3(e_1,e_2)=0,\hfill& \gamma^2_3(e_1,e_3)= -6e_1,\hfill&
\gamma^2_3(e_1,e_4)=-15e_2,\hfill& \gamma^2_3(e_2,e_3)=-9e_2.\cr}$$ We see, in
particular, that $$\langle\gamma^1_2,a^{(2)}\rangle=0+{1\over2}+3-{9\over2}=-
1\ne0,\eqno(4)$$ $$\langle\gamma^2_2,a^{(2)}\rangle=2+{15\over2}-{45\over2}=-
13\ne0.\eqno(5)$$ A non-trivial cycle in $C^{(3)}_2(L_1;L_1^\ast)$ may
be chosen as $$a^{(3)}=e_2^\ast\otimes(e_1\wedge e_4-3e_2\wedge
e_3)$$(again it is a cycle in virtue of a direct calculation, and it is
not homologous to 0 in virtue of the calculation below), and we see that
$$\langle\gamma^1_3,a^{(3)}\rangle=-3+3=0,\eqno(6)$$
$$\langle\gamma^2_3,a^{(3)}\rangle=-15+27=12\ne0.\eqno(7)$$ The
inequalities (5) and (7) show that the deformation $[\, ,\, ]^2_t$ is
singular: one-parameter families of transformations of $L_1$ cannot make
cohomology classes of $\gamma^2_2,\gamma^2_3$ collinear, and no parameter
change can transform this deformation into a deformation with a non-zero
infinitesimal deformation.

On the contrary, the deformation $[\, ,\, ]^1_t$ is non-singular. Indeed,
having applied an appropriate transformation $\varphi_t(g)=g+\lambda(g)t^3$,
we kill $\gamma^1_3$. Then we will have$$\delta\gamma_5^1=-[\gamma^1_1,
\gamma^1_4]-[\gamma^1_2,\gamma^1_3] =0;$$the cocycle $\gamma^1_{(5)}\in
C^2_{(5)}(L_1;L_1)$ is cohomologous to zero (because $H^2_{(5)}(L_1;L_1)=0$),
and we can kill $\gamma^1_{(5)}$ by an appropriate transformation
$\varphi_t(g)=g+\mu(g)t^5$. Proceeding in this way, we kill all
$\gamma^1_k$
with odd $k$, after which we apply the parameter change $u(t)=t^2$. We get the
deformation $$[e_i,e_j]'_u=[e_i,e_j]+\sum_{l\ge1}\gamma_{2l}(e_i,e_j)u^l$$
with $\gamma_2$ being not cohomologous to 0.

This completes the proof of Theorem 4.1.\bigskip

\centerline{\bf 5. Some remarkable cochains of ${\bi L}_{\bf 1}$}\medskip

Let $W$ be the $L_1$-module spanned by $e_j$ with all $j\in\bb Z$ and with the
$L_1$-action $e_i(e_j)=(j-i)e_{i+j}$. It is an extension of the adjoint
representation. Define a cochain $$\mu_k\in C^1_{(k)}(L_1;W),\ k\ge2,$$by the
formula$$\mu_k(e_i)=(-1)^{i+1}{k-1\choose i-2}e_{i-k}.$$ Thus $\mu_k(e_i)=0$
if $i=1$ or $i>k+1$, and $\mu_k(e_2)=-e_{2-k},\, \mu_k(e_{k+1})=(-
1)^ke_1.$\smallskip

{\sc Lemma 5.1}. {\it $\delta\mu_k(e_1,e_i)=0$ for all $i,k$.}\smallskip

{\sc Proof.}$$\eqalign{\delta\mu_k(e_1,e_i)&=\mu_k((i-1)e_{i+1})-
[e_1,\mu_k(e_i)]\cr &=\left((-1)^i(i-1){k-1\choose i-1}-(i-k-1)(-1)^{i+1}{k-
1\choose i-2}\right)e_{i-k+1}\cr &=(-1)^i\left({(i-1)(k-1)!\over(k-i)!(i-1)!}-
{(k-i+1)(k-!)!\over(k-i+1)!(i-2)!}\right)e_{i-k+1}=0.\cr}$$

{\sc Lemma 5.2}. {\it The cochains $\delta\mu_2,\, \delta\mu_3,\, \delta\mu_4$
are cocycles in $C^2(L_1;L_1)\subset C^2(L_1;W)$ not cohomologous to
0.}\smallskip

{\sc Proof}. Since $\delta\mu_k(e_i,e_j)\in{\bb C}e_{i+j-k}$, and if $j>i\ge2$
and $k\le4$ then $i+j-k\ge1$, then $\delta\mu_2,\delta\mu_3,\delta\mu_4\in
C^2(L_1;L_1)$. These cochains are obviously cocycles, and they are not
cohomologous to 0, because no non-zero coboundary in $C^2_{(k)}(L_1;L_1)$ with
$k\ge2$ is annihilated by $e_1$: if $\lambda\in C^1_{(k)}(L_1;L_1), k\ge2$,
then $\lambda(e_i)=a_ie_{i-k},\ a_i=0$ if $i\le k$, and if
$\delta\lambda(e_1,e_i)=\left((i-1)a_{i+1}-(i-k-1)a_i\right)e_{i+1}=0$, then
$(k-1)a_{k+1}=0,ka_{k+2}=0,(k+1)a_{k+3}-a_{k+2}=0,(k+2)a_{k+4}-
2a_{k+3}=0,\dots$, which implies successively that $a_i=0$ for
$i=k+1,k+2,k+3,k+4,\dots$\smallskip

{\sc Corollary.} {\it Any cocycle in $C^2(L_1;L_1)$ is cohomologous to a
linear combination of $\delta\mu_2,\delta\mu_3,\delta\mu_4$}.
\medskip

Define the cochains $\delta_k\in C^2_{(k)}(L_1;L_1)$ by the formula
$$\delta_k(e_i,e_j)=\cases{\delta\mu_k(e_i,e_j)&if $i+j-k\ge1$,\cr
0&otherwise,\cr}$$(in particular, $\delta_k=\delta\mu_k$ if $k=2,3,4$), and
put
$$\delta_{k,l}=[\delta_k,\delta_l]\in C^3_{(k+l)}(L_1;L_1).
$$

{\sc Lemma 5.3}. {\it For $k,l,s,t$ fixed and $N>1+{\rm max}(k,l,k+l-s,k+l-t)$
$$\delta_{k,l}(e_s,e_t,e_N)=a(k,l,s,t)(N+k+l-s-t)e_{N+s+t-k-l}.$$ where
$a(k,l,s,t)$ depends on $k,l,s,t$, but not on $N$.}
\smallskip

{\sc Proof.} Obviously, for $N>k+1$$$\delta\mu_k(e_s,e_N)=(-1)^s{k-1\choose s-
2}(N+k-s)e_{N+s-k}.$$By definition, $\delta_{k,l}(e_s,e_t,e_u)$ is the sum of
6 summands. One of them is
$\delta_k(\delta_l(e_s,e_t),e_N)=\delta\mu_k(\delta_l(e_s,e_t),e_N),$which is
either 0, or$$\eqalign{\delta\mu_k(\delta\mu_l(e_s,e_t),e_N)&=a\mu_k(e_{s+t-
l},e_N)\cr &=ab(N+k+l-s-t)e_{N+s+t-k-l},\cr}$$ 
where $a$ depends only on
$l,s,t$, and $b$ depends only on $k,l,s,t$.  A similar formula is true for
$\delta_k(\delta_l(e_s,e_t),e_N)$. Two more summands are
$$\delta\mu_k(\delta\mu_l(e_t,e_N),e_s)+\delta\mu_l(\delta\mu_k(e_N,e_s),e_t).
$$This sum is $e_{N+s+t-k-l}$ times $$\eqalign{&-(-1)^{s+t}{l-1\choose t-2}{k-
1\choose s-2}(N+l-t)(N+k-l+t-s)\cr &+(-1)^{s+t}{k-1\choose s-2}{l-1\choose t-
2}(N+k-s)(N+l-k+s-t)\cr &\qquad\qquad=(-1)^{s+t}{l-1\choose t-2}{k-1\choose s-
2}(k-t-l+s)(N+k+l-t-s).\cr}$$ 
The last two summands
$$\delta\mu_k(\delta\mu_l(e_N,e_s),e_t)+\delta\mu_l(\delta\mu_k(e_t,e_N),e_s)
$$
are treated similarly.  The proof of Lemma 5.3 is now complete. \smallbreak

{\sc Lemma 5.4}. {\it The quantity $a(k,l,s,t)$ from Lemma
$5.3$ is equal to $0$ if $s+t-k-l>1$.}\smallskip

{\sc Proof.} If $s+t>k+l+1$, then all the summands comprising $a(k,l,s,t)$ are
equal to 0 separately (see the proof of Lemma 5.3).\medskip

Define the cochains $\lambda_{k,l}\in C^2_{(k+l)}(L_1;W)$ by the formula
$$\lambda_{k,l}(e_s,e_t)=a(k,l,s,t)e_{s+t-k-l}.$$ Lemma 5.4 shows that the
cochain $\lambda_{k,l}$ takes values in the subspace of $W$ spanned by $e_i$
with $i\le1$.\smallskip

{\sc Definition.} We say that two cochains $\alpha ,\beta\in C^q(L_1;W)$ are
{\it commensurable}, $\alpha\sim\beta$, if the equality $$\alpha(e_{i_1},\dots
,e_{i_q})=\beta(e_{i_1},\dots ,e_{i_q})$$ holds for all but finitely many sets
$i_1,\dots ,i_q$. In particular,$$\lambda_{k,l}\sim0\ {\rm for\ all\
}k,l.$$\smallskip

{\sc Lemma 5.5}. {\it $\delta\lambda_{k,l}\sim\delta_{k,l}$ for all
$k,l$.}\smallskip

{\sc Proof.} Lemma 5.3 implies that for any $s,t$
$$\delta\lambda_{k,l}(e_s,e_t,e_N)=\delta_{k,l}(e_s,e_t,e_N)$$ when $N$ is
large. It is also true that if $k\le l,\, l+1<s<t<u$, then
$$\delta\lambda_{k,l}(e_s,e_t,e_u)=\delta_{k,l}(e_s,e_t,e_u)=0.$$\smallskip

The following is the most important property of the commensurability
relation.\smallskip

{\sc Lemma 5.6}. {\it If $\alpha ,\beta\in C^q_{(k)}(L_1;W),\
\alpha\sim\beta$, and $\delta\alpha =\delta\beta =0$, then
$\alpha=\beta$.}\smallskip

{\sc Proof.} Let $\gamma=\alpha-\beta$, let $\gamma(e_{i_1},\dots
,e_{i_q})\ne0,$ and let $\gamma(e_{j_1},\dots ,e_{j_q})=0$ if $j_q\ge N$. Take
$j\ge N,j\ne i_1+\dots+i_q-k$; then $$\delta\gamma(e_{i_1},\dots
,e_{i_q},e_j)=\pm[e_j,\gamma(e_{i_1},\dots ,e_{i_q})]\ne0$$ which contradicts
to $\delta\gamma=0$.\medskip

With the exception of $\lambda_{2,2}\in C^2_{(4)}(L_1;L_1)$, the cochains
$\lambda_{k,l}$ are not contained in $C^2(L_1;L_1)$. But in some cases we can
force them into $C^2(L_1;L_1)$ by adding an appropriate multiple of
$\delta\mu_{k+l}$ (which would not affect the coboundaries).\smallskip

{\sc Lemma 5.7}.$$\eqalign{\lambda_{2,2}&\in C^2_{(4)}(L_1;L_1)\cr
13\lambda_{2,3}-\delta\mu_5&\in C^2_{(5)}(L_1;L_1)\cr
7\lambda_{2,4}+2\delta\mu_6&\in C^2_{(6)}(L_1;L_1)\cr 21\lambda_{3,3}-
10\delta\mu_6&\in C^2_{(6)}(L_1;L_1)\cr 65\lambda_{2,5}+119\lambda_{3,4}-
2\delta\mu_7&\in C^2_{(7)}(L_1;L_1)\cr 151\lambda_{2,6}-
105\lambda_{3,5}+60\delta\mu_8&\in C^2_{(8)}(L_1;L_1)\cr
20\lambda_{2,6}+7\lambda_{4,4}+6\delta\mu_8&\in C^2_{(8)}(L_1;L_1)\cr}$$

{\sc Proof:} direct calculation. Each of the cochains
$\lambda_{k,l},\delta\mu_m$ has finitely many values outside $L_1$, and it is
very easy to compute them all. For example,$$\lambda_{2,3}(e_2,e_3)=2e_0,\
\delta\mu_5(e_2,e_3)=26e_0,$$all other values of these cochains are in $L_1$.
Hence $13\lambda_{2,3}-\delta\mu_5\in C^2(L_1;L_1)$. Similarly
$$\matrix{\lambda_{2,4}(e_2,e_3)=-12e_{-1},\hfill
&\lambda_{3,3}(e_2,e_3)=20e_{-1},\hfill&\delta\mu_6(e_2,e_3)=42e_{-
1},\hfill\cr \lambda_{2,4}(e_2,e_4)=12e_0,\hfill&\lambda_{3,3}(e_2,e_4)=-
20e_0,\hfill&\delta\mu_6(e_2,e_4)=-42e_0,\hfill\cr}$$all other values of 
these cochains are in $L_1$. Hence $7\lambda_{2,4}+2\delta\mu_6\in 
C^2(L_1;L_1),\ 21\lambda_{3,3}-10\delta\mu_6\in C^2(L_1;L_1)$. In the 
remaining cases the computations are similar (though longer).\smallskip 

{\sc Lemma 5.8}. {\it If a linear combination of $\delta_{2,5}$ and
$\delta_{3,4}$ is commensurable with a coboundary, then this linear
combination is a multiple of $65\delta_{2,5}+119\delta_{3,4}$. If a linear
combination of $\delta_{2,6}, \delta_{3,5}$ and $\delta_{4,4}$ is
commensurable with a coboundary, then this linear combination is a linear
combination of $151\delta_{2,6}-105\delta_{3,5}$ and
$20\delta_{2,6}+7\delta_{4,4}$. In other words, if
$a\delta_{2,6}+b\delta_{3,5}+c\delta_{4.4}$ is commensurable with a
coboundary, then $105a+151b-300c=0$}.\smallskip

{\sc Proof.} Lemma 5.7 implies that $65\delta_{2,5}+119\delta_{3,4},
151\delta_{2,6}-105\delta_{3,5}$ and $20\delta_{2,6}+7\delta_{4,4}$ are
commensurable with coboundaries, and Lemma 5.6 and Theorem 3.2 imply that
$\delta_{3,4}$ and $\delta_{4,4}$ are not commensurable with coboundaries.
The Lemma follows.\smallskip

In conclusion we remark, that adding linear combinations of $\lambda$'s to
$\delta\mu$'s does not affect essentially the Lie products. For
example,$$[13\lambda_{2,3}-\delta\mu_5,7\lambda_{2,4}+2\delta\mu_6]\sim-
2\delta_{5,6},$$and so on.\bigskip

\centerline{\bf 6. Proof of Theorem 1.1}\medskip

We consider formal one-parameter deformations $$[g,h]_t=[g,h]+\sum_{k\ge1}
\alpha_k(g,h)t^k$$of $L_1$; recall that the relation (1) from Section 2 should
hold for all the cochains $\alpha_k\in C^2(L_1;L_1)$.

The general theory of Section 2 combined with Theorem 3.1 and the Corollary to
Lemma 5.2 imply that the following construction gives all possible formal
deformations of $L_1$ up to the equivalence described in Section 2.

Suppose that $\alpha_1,\dots ,\alpha_{k-1}$ has already been defined. Fix an
arbitrary cochain $\alpha^0_k\in C^2(L_1;L_1)$ with
$$
\delta\alpha^0_k=-
{1\over2}\sum_{i=1}^{k-1}[\alpha_i,\alpha_{k-i}],\eqno (8)
$$
then choose three complex numbers $c_{k2},c_{k3},c_{k4}$, and put
inductively
$$\alpha_k=\alpha_k^0+c_{k2}\delta\mu_2+c_{k3}\delta\mu_3+c_{k4}\delta\mu_4.$$
Notice that we do not vary the cochains $\alpha^0_k$, their choice is
arbitrary, but after having been chosen, they are fixed. On the contrary,
the numbers $c_{ij}$ are varied, but the existence of $\alpha^0_k$
imposes a condition on all $c_{ij}$ with $i<k$. We will always choose
$\alpha^0_k$ under the condition: if the right hand side of (8) lies in
$\oplus_{q\le m}C^3_{(q)}(L_1;L_1)$, then $\alpha^0_k\in\oplus_{q\le
m}C^2_{(q)}(L_1;L_1)$.

Below we use the notation: $\oplus_{q\le m}C^r_{(q)}(L_1;L_1)=C^r_{(\le
m)}(L_1;L_1)$.

Our goal is to show that the conditions imposed on $c_{ij}$ by the solvability
of the equations (8) leave only a few possibilities for the
deformation.\smallskip

{\sc Lemma 6.1.} $c_{13}=c_{14}=0.$\smallskip

This follows from Theorem 3.2. The argumentation is well known (see [Fu]), but
we give it here for the completeness sake. The degree 8 component of
$[\alpha_1,\alpha_1]$ is $c_{14}^2[\delta\mu_4,\delta\mu_4]$, and since it
belongs to Im$\, \delta$ and $[c,c]\ne0$, we must have $c_{14}=0$. Let $\beta$
be the component of degree 6 of $\alpha_2$; then $\delta\beta=-
\displaystyle{1\strut\over\strut2}c_{13}^2[\delta\mu_3,\delta\mu_3]$. Since
the degree 7 component of $\alpha_2$ is 0, the degree 9 component of
$[\alpha_1,\alpha_2]$ is $c_{13}[\delta\mu_3,\beta]$, and since $\langle
b,b,b\rangle\ne0$, the latter belongs to Im$\, \delta$ only if
$c_{13}=0$.\smallskip

It remains to consider two cases: $c_{12}\ne0$ and $c_{12}=0$. The parameter
change $u=c_{12}t$ would reduce the case $c_{12}\ne0$ to the case $c_{12}=1$.
Consider this case.\smallskip

{\sc Lemma 6.2.} {\it Let $\alpha_1=\delta\mu_2$. Then $\alpha_k\in C^2_{(\le
2k)}(L_1;L_1).$}\smallskip

{\sc Proof.} Let $k\ge2$ and suppose that $\alpha_l\in C^2_{(\le
2l)}(L_1;L_1)$ for $l<k$. Then the right hand side of (8) lies in $C^3_{(\le
2k)}(L_1;L_1)$, and hence $\alpha_k^0\in C^2_{(\le 2k)}(L_1;L_1)$. But since
$k\ge2$, then $\delta\mu_2,\delta\mu_3,\delta\mu_4\in C^2_{(\le
2k)}(L_1;L_1)$. Thus $\alpha_k\in C^2_{(\le 2k)}(L_1;L_1).$\medskip

Consider now the ``homogeneous" case.\smallskip

{\sc Lemma 6.3.} {\it There exist $($up to an equivalence$)$ at most two
deformations with $\alpha_1=\delta\mu_2$ and $\alpha_k\in
C^2_{(2k)}(L_1;L_1).$}\smallskip

{\sc Proof.} First, $\displaystyle{\delta\alpha_2=-
{1\strut\over\strut2}[\alpha_1,\alpha_1]=-{1\over2}[\delta\mu_2,\delta\mu_2]=-
{1\strut\over\strut2}\delta_{2,2}}.$ Since $\delta_{2,2}=\delta\lambda_{2,2}$
(Lemmas 5.5 and 5.6), we can take $\alpha_2^0=-
\displaystyle{1\strut\over\strut2}\lambda_{2,2}$, and then we have $\alpha_2=-
\displaystyle{1\strut\over\strut2}\lambda_{2,2}+x\delta\mu_4$. Notice, that in
the whole construction $x$ is the only constant to choose: if $k\ge3$, then
deg$\, \alpha_k>4$, and $\alpha_k$ does not involve either of
$\delta\mu_2,\delta\mu_3,\delta\mu_4$. Next we have $\delta\alpha_3=-
[\alpha_1,\alpha_2]=-\left[\delta\mu_2,-
\displaystyle{1\strut\over\strut2}\lambda_{2,2}+x\delta\mu_4\right]\sim-
x\delta_{2,4}$, and we may put $\alpha_3=-x\lambda_{2,4}-
\displaystyle{2\strut\over\strut7}x\delta\mu_6$ (Lemmas 5.7 and 5.6). Further
we have $$\eqalign{\delta\alpha_4&=-[\alpha_1,\alpha_3]-
{1\over2}[\alpha_2,\alpha_2]\cr &=-\left[\delta\mu_2,-x\lambda_{2,4}-
{2\over7}x\delta\mu_6\right]-{1\over2}\left[-
{1\over2}\lambda_{2,2}+x\delta\mu_4,-
{1\over2}\lambda_{2,2}+x\delta\mu_4\right]\cr &\sim {2\over7}x\delta_{2,6}-
{1\over2}x^2\delta_{4,4}.}$$ According to Lemma 5.8 the latter should be a
linear combination of $151\delta_{2,6}-105\delta_{3,5}$ and
$20\delta_{2,6}+7\delta_{4,4}$, that is $$105\cdot{2\over7}x-
300\cdot{1\over2}x^2=210x(5x-1)=0.$$ Hence $x=0$ or 
$-\displaystyle{1\strut\over\strut5}$.  This completes the proof of Lemma 6.3.
\smallskip 

Actually, the two deformations do exist: the deformation $[\, ,\, ]^3_t$ and
the deformation constructed in Section 4 from $[\, ,\, ]^1_t$. We do not check
it now: it will be easier to see later.\smallskip

{\sc Lemma 6.4}. {\it Up to a parameter change, there are no deformations with
$\alpha_1=\delta\mu_2$, other than those in Lemma} 6.3.\smallskip

{\sc Proof.} First we notice that using an appropriate parameter change, we
may kill $\delta\mu_2$ in
$\alpha_k=\alpha_k^0+c_{k2}\delta\mu_2+c_{k3}\delta\mu_3+c_{k4}\delta\mu_4$
for all $k\ge2$. Hence we can write:$$\eqalign{\alpha_1&=\delta\mu_2,\cr
\alpha_2&=-{1\over2}\lambda_{2,2}+y_2\delta\mu_3+x_2\delta\mu_4,\cr
\delta\alpha_3&\sim -y_2\delta_{2,3}-x_2\delta_{2,4}\cr \alpha_3&=-
y_2\lambda_{2,3}+{1\over13}y_2\delta\mu_5-x_2\lambda_{2,4}-
{2\over7}x_2\delta\mu_6+y_3\delta\mu_3+x_3\delta\mu_4,\cr \delta\alpha_4&\sim
-{1\over2}y_2^2\delta_{3,3}-x_2y_2\delta_{3,4}-{1\over2}x_2^2\delta_{4,4}-
{2\over13}y_2\delta_{2,5}+{2\over7}x_2\delta_{2,6}-y_3\delta_{2,3}-
x_3\delta_{2,4}.\cr}$$ As in the previous proof, $\displaystyle{-
{1\strut\over\strut2}x_2^2\delta_{4,4}+{2\over7}x_2\delta_{2,6}}$ being in
Im$\, \delta$ implies (by Lemma 5.8) that either $x_2=0$ or $x_2=-
\displaystyle{1\strut\over\strut5}$. Consider these two cases.\smallskip

{\it Case} 1: $x_2=0$. In this case $$\delta\alpha_4\sim -
{1\over2}y_2^2\delta_{3,3}-{2\over13}y_2\delta_{2,5}-y_3\delta_{2,3}-
x_3\delta_{2,4},$$ and since $\displaystyle{-
{2\strut\over\strut13}y_2\delta_{2,5}}\in{\rm Im}\, \delta$, Lemma 5.8 implies
that $y_2=0$. Thus $\alpha_1=\delta\mu_2,
\alpha_2=\displaystyle{1\strut\over\strut2}\lambda_{2,2},\alpha_3=
y_3\delta\mu_3 +x_3\delta\mu_4$, and hence $$\eqalign{\alpha_4&=-
y_3\lambda_{2,3}+{1\over13}y_3\delta\mu_5-x_3\lambda_{2,4}-
{2\over7}x_3\delta\mu_6+y_4\delta\mu_3+x_4\delta\mu_4,\cr \delta\alpha_5&\sim
-{1\over13}y_3\delta_{2,5}+{2\over7}x_3\delta_{2,6}-y_4\delta_{2,3}-
x_4\delta_{2,4},\cr}$$ which implies, in virtue of Lemma 5.8, that
$y_3=0,x_3=0$. Similarly, for $k>4$ we have inductively$$\eqalign{\alpha_k&=-
y_{k-2}\lambda_{2,3}+{1\over13}y_{k-1}\delta\mu_5-x_{k-1}\lambda_{2,4}-
{2\over7}x_{k-1}\delta\mu_6+y_k\delta\mu_3+x_k\delta\mu_4,\cr
\delta\alpha_{k+1}&\sim -{1\over13}y_{k-1}\delta_{2,5}+{2\over7}x_{k-
2}\delta_{2,6}-y_k\delta_{2,3}-x_k\delta_{2,4},\cr}$$ which implies that
$y_{k-2}=0,x_{k-2}=0$.\smallskip

{\it Case} 2: $x_2=-\displaystyle{1\strut\over\strut5}$. In the calculations
below we skip for all $\alpha$'s the terms of degree $\le8$ and for all
$\delta\alpha$'s the terms of degree $\le9$. We
have:$$\eqalign{\alpha_1&=\delta\mu_2,\cr \alpha_2&=-{1\over2}\lambda_{2,2}-
{1\over5}\delta\mu_4,\cr \alpha_3&=-
{1\over5}\lambda_{2,4}+{2\over35}\delta\mu_6+y_3\delta\mu_3+x_3\delta\mu_4,\cr
\delta\alpha_4&\sim -{1\over50}\delta_{4,4}+{2\over35}\delta_{2,6}-
y_3\delta_{2,3}-x_3\delta_{2,4},\cr \alpha_4&=\dots -
y_3\lambda_{2,3}+{1\over13}y_3\delta\mu_5-x_3\lambda_{2,4}-
{2\over7}x_3\delta\mu_6+y_4\delta\mu_3+x_4\delta\mu_4,\cr
\delta\alpha_5&\sim\dots -{1\over13}y_3\delta_{2,5}+{2\over7}x_3\delta_{2,6}-
y_4\delta_{2,3}-
x_4\delta_{2,4}+{1\over5}y_3\delta_{3,4}+{1\over5}x_3\delta_{4,4}.\cr}$$ The
latter implies that $x_3=0,y_3=0$: otherwise $-\displaystyle
{1\strut\over\strut13}y_3\delta_{2,5}$ and
$\displaystyle{{2\strut\over\strut7}x_3\delta_{2,6}
+{1\over5}x_3\delta_{4,4}}$ belonging to Im$\, \delta$ contradicts to Lemma
5.8. Further, for $k>4$ we have inductively $$\eqalign{\alpha_k&=\dots -y_{k-
1}\lambda_{2,3}+{1\over13}y_{k-1}\delta\mu_5-x_{k-1}\lambda_{2,4}-
{2\over7}x_{k-1}\delta\mu_6+y_k\delta\mu_3+x_k\delta\mu_4,\cr
\delta\alpha_{k+1}&\sim\dots -{1\over13}y_{k-1}\delta_{2,5}+{2\over7}x_{k-
1}\delta_{2,6}-y_k\delta_{2,3}-x_k\delta_{2,4}+{1\over10}y_{k-
1}\delta_{3,4}+{1\over10}x_{k-1}\delta_{4,4},\cr}$$which implies $x_{k-
1}=0,y_{k-1}=0.$

Thus in both cases we have no deformations different from those of Lemma
6.3. This completes the proof of Lemma 6.4.\smallskip

The case $\alpha_1\ne0$ is over; consider the case $\alpha_1=0$. In this case
we are interested only in {\it singular} deformations (See Section 4).

First of all, $\delta\alpha_2=\delta\alpha_3=0$, and we can assume that
$$\eqalign{\alpha_2&=z_2\delta\mu_2+y_2\delta\mu_3+x_2\delta\mu_4,\cr
\alpha_3&=z_3\delta\mu_2+y_3\delta\mu_3+x_3\delta\mu_4.\cr}$$ The degree 8
component of $\delta\alpha_4$ should be 
$-\displaystyle{1\strut\over\strut2}x_2^2[\delta\mu_4,\delta\mu_4]$, 
which is possible only if $x_2=0$ (Theorem 3.2). Now, the degree 9 
component of $\delta\alpha_6$ will be a cocycle representing the 
Massey-Lie cube of the cohomology class of $\delta\mu_3$ times some 
non-zero factor times $y_2^2$, which implies that $y_2=0$ (again Theorem 
3.2). After it, applying Theorem 3.2 to the degree 8 component of 
$\alpha_6$ we find that $x_3=0$. 

Thus $\alpha_2=z_2\delta\mu_2, \alpha_3=z_3\delta\mu_2+y_3\delta\mu_3$. As
before, two essentially different cases are $z_2=0$ and $z_2\ne0$. Suppose
that $z_2\ne0$. Then using an appropriate parameter change, we can make
$z_2=1, z_3=0$. Thus we have $\alpha_2=\delta\mu_2, \alpha_3=y\delta\mu_3$
(here $y=y_3$); also we have $\alpha_k\in C^2_{(\le k)}(L_1;L_1)$, which is
similar to Lemma 6.2. As before, we begin with the ``homogeneous"
case.\smallskip

{\sc Lemma 6.5.} {\it There exist $($up to an equivalence$)$ at most one
singular deformation with $\alpha_1=0, \alpha_2=\delta\mu_2$ and $\alpha_k\in
C^2_{(k)}(L_1;L_1)$}.\smallskip

{\sc Proof.} We have $\alpha_1=0, \alpha_2=\delta\mu_2,
\alpha_3=y\delta\mu_3,$$$\delta\alpha_4=-{1\over2}\delta_{2,2}, \alpha_4=-
{1\over 2}\lambda_{2,2}+x\delta\mu_4\sim x\delta\mu_4,$$and $y$ and $x$ are
the only parameters of which our deformation can depend. Next we have
$$\delta\alpha_5=-y\delta_{2,3},\alpha_5=-
\lambda_{2,3}+{1\over13}\delta\mu_5\sim{1\over13}y\delta\mu_5,$$
$$\delta\alpha_6\sim -x\delta_{2,4}-{1\over2}y^2\delta_{3,3},\alpha_6\sim
\left({5\over21}y^2-{2\over7}x\right)\delta\mu_6,$$$$\delta\alpha_7\sim -
{1\over13}y\delta_{2,5}-xy\delta_{3,4}.$$According to Lemma 5.8, the latter
implies the relation$${1\over13}y:xy=65:119,$$which may mean either $y=0$ or
$x=\displaystyle{119\strut\over\strut{13\cdot65}}$. But the first would imply
$\alpha_3=0$ and successively $\alpha_5=0,\alpha_7=0,\dots$, which means that
the deformation is actually non-singular. This leaves us with the second
possibility. Next we have:$$\delta\alpha_8\sim\left(-
{5\over21}y^2+{2\over7}x\right)\delta_{2,6}-{1\over13}y^2\delta_{3,5}-
{1\over2}x^2\delta_{3,4},$$and we again apply Lemma 5.8 to get the second
relation:$$105\left({5\over21}y^2-{2\over7}x\right)+151\cdot {1\over13}y^2-
300\cdot{1\over2}x^2=0,$$ which determines $y$ up to the sign:
$y^2=\displaystyle{2\strut\cdot6^3\over\strut13^3}$. The sign of $y$ is
irrelevant, for it may be changed by the parameter change $t\mapsto -t$.
This completes the proof of Lemma 6.5. \smallskip

{\sc Lemma 6.6}. {\it There are no more one-parameter formal deformations of
$L_1$. In other words, any one-parameter formal deformation of $L_1$ is
reduced by an equivalence and a parameter change to one of the deformations of
Lemmas $6.3$ and $6.5$.}\smallskip

{\sc Proof.} First we consider the case of a non-homogeneous deformation with
$\alpha_1=0, \alpha_2=\delta\mu_2, \alpha_3=y\delta\mu_3,
\alpha_4=y_4\delta\mu_3+x\delta\mu_4$. In this case
$\alpha_5\sim\displaystyle{1\strut\over\strut13}y\delta\mu_5+x_5\delta\mu_4
+y_5\delta\mu_3$, and the following formulas give the components of the degree
less by one than the maximal one:$$\eqalign{\alpha_6&\sim\dots-
y_4\delta_{2,3}\dots,\cr \delta\alpha_6&\sim\dots+{1\over13}\delta\mu_5\dots,
\cr \delta\alpha_7&\sim\dots-x_5\delta_{2,4}-yy_4\delta_{3,3}\dots,\cr
\alpha_7&\sim\dots\left(-{2\over7}x_5+{10\over21}yy_4\right)\delta\mu_6\dots,
\cr \delta\alpha_8&\sim\dots-{1\over13}y_4\delta_{2,5}-xy_4\delta_{3,4}-
x_5y\delta_{3,4}\dots,\cr \delta\alpha_9&\sim\dots\left({2\over7}x_5-
{10\over21}yy_4\right)\delta_{2,6}-{2\over13}yy_4\delta_{3,5}-
xx_5\delta_{4,4}.\cr}$$The last two equalities, combined with Lemma 5.8 give
two linear relations between $y_4$ and $x_5$: $$\eqalign{\left({119\over13}-
65x\right)y_4+65yx_5=0,\cr {8\cdot119\over13}yy_4-(30+300x)x_5=0.\cr}$$The
calculations made above leave us with only three possibilities for $x$ and
$y$:$$y=x=0;\ y=0,x=-{1\over5};\ y={12\over13}\sqrt{3\over13},x={119\over13
\cdot65};$$in each of these cases the determinant of the above system is not
equal to zero, and hence $y_4=x_5=0$. Then we consider the components of
degree 7 and 8 for $\delta\alpha_9$ and $\delta\alpha_{10}$, and we obtain
precisely the same system for $y_5$ and $x_6$; hence they are also equal to
zero. Proceeding in this way we prove that our deformation is actually
homogeneous and is covered by Lemma 6.5.

And the last case is $\alpha_1=0,\alpha_2=0$. Suppose that $\alpha_3=0,\dots,
\alpha_{p-1}=0,\alpha_p\ne0$ and that the deformation is not equivalent to a
deformation with $\alpha_1=0,\dots\alpha_p=0$. In this case $\delta\alpha_i=0$
and $\alpha_i-z_i\delta\mu_2+y_i\delta\mu_3+x_i\delta\mu_4$ for $p\le i<2p-
1$. If $x_i=0$ for $i<p+r$ for some $r,\, 0\le r<p$, then the degree 8
component first appears in $\delta\alpha_{2p+2r}$ and is equal to 
$-\displaystyle{1\strut\over\strut2}x_{p+r}^2\delta_{4,4}$, which implies 
that $x_{p+r}=0$. Hence all $x_i=0$, and 
$\alpha_i=z_i\delta\mu_2+y_i\delta\mu_3$ for $p\le i<2p-1$. Now, the 
degree 9 component first appear in $\delta\alpha_{3p}$ and is 
$\sim{5\strut\over\strut21}y_p^3\delta_{3,6}$, which implies in virtue of 
Lemma 5.8, that $y_p=0$. This equality shifts the first appearance of the 
degree 9 component to $\delta\alpha_{3p+3}$, and Lemma 5.8 yields 
$y_{p+1}=0$. Proceeding in the same way, we prove that $y_i=0$ for 
$i<\displaystyle{3\strut\over\strut2}p$. A parameter change makes 
$\alpha_p=\delta\mu_2$ and kills $\delta\mu_2$ term in all $\alpha_i$'s 
with $i>p$; in particular, it makes $\alpha_i=0$ for 
$p<i<\displaystyle{3\strut\over\strut2}p$. 

In the remaining part of the proof the cases of even and odd $p$ are slightly
different. If $p=2q$, then we have $\alpha_{2q}=\delta\mu_2,\,
\alpha_{3q}=y\delta\mu_3,\, \alpha_i=y_i\delta\mu_3\ {\rm if}\ 3q<i<4q,\,
\alpha_{4q}=y_{4q}\delta\mu_3+x\delta\mu_4,\,
y_{4q+1}=y_{4q+1}\delta\mu_3+x_{4q+1}\delta\mu_4$. Then we consider the
degree 7 component of $\delta\alpha_{7q}$ and the degree 8 component of
$\delta\alpha_{8q}$, and, using Lemma 5.8, obtain the same system of equations
for $y$ and $x$ as in the proof of Lemma 6.5; we already know all the
solutions of this system. Next we consider the degree 7 component of
$\delta\alpha_{7q+1}$ and the degree 8 component of $\delta\alpha_{8q+1}$, and
obtain for $y_{3q+1}$ and $x_{4q+1}$ the same equations as for $y_4$ and $x_5$
in the first part of this proof; these equations imply that
$y_{3q+1}=x_{4q+1}=0$. Proceeding in the same way, we annihilate all terms
with $\delta\mu_3$ and $\delta\mu_4$, with the exception of $y\delta\mu_3$ in
$\alpha_{3q}$ and $x\delta\mu_4$ in $\alpha_{4q}$, and finally we get
$\alpha_i=0$ if $i$ is not divisible by $q$. The parameter change $t^q\mapsto
t$ makes our deformation the deformation of Lemma 6.5 if $y\ne0$ and one of
the deformations of Lemma 6.3 if $y=0$.

If $p$ is odd, then we proceed as above and get
$\displaystyle{\alpha_p=\delta\mu_2,\, \alpha_i=0\ {\rm if}\ p+1\le
i<{3p+1\strut\over\strut2}},$ $\displaystyle{\alpha_i=y_i\delta\mu_3\ \ {\rm
if}\ \ {3p+1\strut\over\strut2}\le i<2p,}$ $\displaystyle{
\alpha_{2p}=y_{2p}\delta\mu_3+x\delta\mu_4,}$ $\displaystyle{
\alpha_i=y_i\delta\mu_3+x_i\delta\mu_4\ \ {\rm if}\ \ 2p<i<}$\break
$\displaystyle{{5p+1\strut\over\strut2},\ \
\alpha_i={1\strut\over\strut13}y_{i-p}\delta\mu_5+y_i\delta\mu_3+
x_i\delta\mu_4\ \ {\rm if}\ \ {1\strut\over\strut13}\le i<3p,}$ 
$\displaystyle{\alpha_{3p}=-{2\strut\over\strut7}x\delta\mu_6+ 
{1\strut\over\strut13}y_{2p}\delta\mu_5+}$
\break 
$\displaystyle{y_i\delta\mu_3 +x_i\delta\mu_4,\, \alpha_{3p+1}=-
{2\strut\over\strut7}x_{2p+1}\delta\mu_6+ 
{1\strut\over\strut13}y_{2p+1}\delta\mu_5}$. The degree 8 component of 
$\delta\alpha_{4p}$ is 
$\sim\displaystyle{{2\strut\over\strut7}x\delta_{2,6} -
{1\strut\over\strut2}x^2\delta_{4,4}}$, and Lemma 5.8 implies $x=0$ or 
$-\displaystyle{1\strut\over\strut5}$ (compare the proof of Lemma 6.3). 
Further the degree 7 component of $\delta\alpha_{7p+1\over2}$ is 
$\sim\displaystyle{-{1\strut\over\strut13}y_{3p+1\over2}\delta_{2,5}-
{1\strut\over\strut2}xy_{3p+1\over2}\delta_{3,4}}$. For the both values 
of $x$ Lemma 5.8 implies $y_{3p+1\over2}=0$. Further we consider the 
degree 8 component of $\delta\alpha_{8p+1}$, the degree 7 component of 
$\delta\alpha_{7p+3\over2}$, the degree 8 component of 
$\delta\alpha_{8p+2}$, and so on, and Lemma 5.8. implies that 
$x_{2p+1}=0,\,  y_{3p+3\over2}=0,\, x_{2p+2}=0,$ and so on. Finally we 
get $\alpha_i=0$ if $i$ is not divisible by $p$. The parameter change 
$t^p\mapsto t$ makes our deformation the deformation of Lemma 6.3. 

This completes the proof of Lemma 6.6.\smallskip

Lemmas 6.3 -- 6.6 show that there are at most three different formal one
parameter deformations of $L_1$, and we already know, that three such
deformations exist. This ends the proof of Theorem 1.1.

\vskip 30pt

\centerline{BIBLIOGRAPHY}\bigskip

\item{[FeFu]} Feigin B., Fuchs D., Homology of the Lie algebras of vector
fields in the line. Funct. Anal. and appl., 14:3 (1980), 45--60

\item{[Fi1]} Fialowski, A., Deformations of the Lie algebra of vector fields
on the line, Uspekhi Mat. Nauk, 38 (1983), pp.201--202; English
translation: Russian Math. Surveys, 38 (1983), No.1, 185-186

\item{[Fi2]} Fialowski, A., An example of formal deformations of Lie
algebras, NATO Conference on Deformation Theory of Algebras and
Applications, Il Ciocco, Italy, 1986. Proceedings, Kluwer, Dordrecht,
1988, pp. 375-401

\item{[Fi3]} Fialowski A., Deformations of Lie algebras, Mat. Sb. (N.S.)
127 (169), (1985), 476-482;  English translation: Math. USSR-Sbornik,
55:2 (1986), 467--473

\item{[FiFu]} Fialowski A., Fuchs D., Construction of miniversal 
deformations of Lie algebras (preprint, 1996).

\item{[Fu]} Fuchs D., Cohomology of infinite-dimensional Lie algebras.
Consultants Bureau, NY, Lond., 1986

\item{[FL]} Fuchs D., Lang L., Massey Products and Deformations, 
submitted for publication 

\item{[Ko]} Kontsevich M., Lecture Notes, Berkeley, 1994.

\item{[R]} Retakh V., The Massey operations in Lie superalgebras and 
deformations of complex-analytic algebras, Funct. Anal. and Appl. 11 
(1977), 319-321

\bigskip
\parindent=0pt
\obeylines {\sc E-mail}
\quad Alice Fialowski: {\tt fialowsk@cs.elte.hu}
\quad Dmitry Fuchs: {\tt fuchs@math.ucdavis.edu}
\bye